\documentclass[11pt]{article}
\usepackage{latexsym}
\usepackage{epsfig}
\usepackage{graphicx}
\usepackage{amsfonts}
\usepackage{amsmath}
\usepackage{amsthm}
\font\bbfont=msbm10 at 11pt
\newtheorem{theorem}{Theorem}[section]
\newtheorem{lemma}[theorem]{Lemma}
\newtheorem{proposition}[theorem]{Proposition}

\numberwithin{equation}{section}

\def\R{\mbox{\bbfont R}}
\def\e{\varepsilon}
\def\D{\Delta}
\def\spanr{\mathrm{span}}
\def\f{\bar f}

\def\P{P}

\begin{document}

\title {Vector Spaces Spanned by the Angle Sums of Polytopes}
\author{Kristin A. Camenga \\ Cornell University, Ithaca, NY \\ kacam@math.cornell.edu\footnote{This work was partially supported by NSF grants 0100323 and 9983660.}}
\maketitle

\begin{abstract} This paper describe the spaces spanned by the angle sums of certain classes of polytopes, as recorded in the $\alpha_{}$-vector.  Families of polytopes are constructed whose angle sums span the spaces of polytopes defined by the Gram and Perles equations, analogs of the Euler and Dehn-Sommerville equations. This shows that the dimension of the affine span of the space of angle sums of simplices is $\left\lfloor \frac{d-1}{2}\right\rfloor$, and that of the combined angle sums and face numbers of simplicial polytopes and general polytopes are $d-1$ and $2d-3$, respectively. A tool used in proving these results is the $\gamma$-vector, an angle analog to the $h$-vector. 
\end{abstract}

\section{Introduction}

One of the motivating questions in the combinatorial study of polytopes is whether a given set of combinatorial data arises from a given class of polytopes. This has been studied in depth using the $f$-vector, which counts the number of faces of a polytope of each dimension.  We will study polytopes by considering their angle sums, which quantify a geometric aspect of polytopes.  We will also introduce a method to construct polytopes while controlling the angle sums.

Let $P$ be any $d$-polytope.  For any face $F$ of $P$, we consider a $d$-dimensional ball centered at an interior point of $F$, small enough that it only intersects faces which contain $F$. The \emph{interior angle} at $F$ in $P$, denoted by $\alpha(F,P)$, is the fraction of this ball that is contained in $P$. Therefore, for a $(d-1)$-dimensional facet $F$ of $P$, $\alpha(F,P) = \frac{1}{2}$ and $\alpha(P, P) = 1$.  The \emph{angle sums} of $P$ are defined for $0\leq i\leq
d$ as

\begin{equation*}
\alpha_{i}(P)= \sum_{i-\text{faces } F\subseteq P} \alpha(F,P) .
\end{equation*}

We write $f_{i}(P)$ for the number of $i$-faces of P and define the following for a $d$-polytope $P$:

\hspace{.5in}$\alpha{}$-vector: $(\alpha_{0}(P), \alpha_{1}(P), \ldots , \alpha_{d}(P))$,

\hspace{.5in}$f$-vector: $(f_{0}(P), f_{1}(P), \ldots , f_{d}(P))$,

\hspace{.5in}$\alpha{}$-$f$-vector: $(\alpha_{0}(P), \alpha_{1}(P), \ldots ,
\alpha_{d}(P),f_{0}(P), f_{1}(P), \ldots , f_{d}(P))$.

By convention, we write $f_{-1}(P) = 1$ and $\alpha_{-1}(P) = 0$.  It is well-known that there are equations on the $f$-vector.  For general polytopes, it is known that the only linear relation on the $f$-vector is the Euler relation \cite{Gr}:
\begin{equation*}
\sum_{i=0}^{d} (-1)^{i}f_{i}(P) = 1 \text{ for any }d\text{-polytope }P.
\end{equation*}

\noindent The only relations on the $f$-vectors of simplicial polytopes are the Dehn-Sommerville relations \cite{Gr}:  for any simplicial polytope $P$ and $-1 \leq k \leq d-2$,
\begin{equation*}
\sum_{j=k}^{d-1}(-1)^{j}{j+1 \choose k+1} f_{j}(P) = (-1)^{d-1} f_{k}(P).
\end{equation*}

These relations are frequently written using the \emph{$h$-vector}, a linear transformation of the $f$-vector.  The $h$-vector is defined on a simplicial polytope $\P$ as
$h(\P) = \left(h_0(\P), h_1(\P), \ldots, h_d(\P)\right)$, where
\begin{equation*}\label{h-vector}
h_i(\P) = \sum_{j=0}^{i}(-1)^{i-j}{d-j \choose d-i}f_{j-1}(\P).
\end{equation*}
  Using the $h$-vector, the Dehn-Sommerville relations have a more symmetric reformulation \cite{Ziegler}: for any simplicial polytope $P$ and $i = 0, \ldots, \left\lfloor\frac{d}{2}\right\rfloor$, $h_i(\P) = h_{d-i}(\P)$.
Since this transformation is invertible, the linear independence of a set of $f$-vectors is equivalent to the linear independence of the corresponding set of $h$-vectors.

It is also known that there are relations on the $\alpha$-vector.  The Gram relation is an analog of the Euler relation \cite{Gr}:
\begin{equation*}
\sum_{i=0}^{d} (-1)^{i}\alpha_{i}(P) = 0 \text{ for any }d\text{-polytope }P.
\end{equation*}
\noindent H\"ohn \cite{Hohn} first showed that these equations are the only linear equations on the $\alpha{}$-vector of general polytopes. Perles proved an analog of the Dehn-Sommerville relations \cite{Gr,PS}:  for any simplicial polytope $P$ and $0 \leq k \leq d-1$, we have
\begin{equation*}
\sum_{j=k}^{d-1} (-1)^{j} {j+1\choose k+1} \alpha_{j}(P) =
(-1)^{d}(\alpha_{k}(P) - f_{k}(P)).
\end{equation*}

In the next section, we will define two constructions on polytopes and in the fourth section we will use these constructions to make sets of polytopes whose $\alpha$- or $\alpha$-$f$-vectors span the spaces defined by the relations above on various classes of polytopes.  The proofs in section 4 will be aided by the $\gamma$-vector, an angle analog of the $h$-vector, which is defined and examined in section 3.

\section{Construction of Polytopes}

We will define two constructions, the pyramid and prism operations, that create polytopes with varying angle sums. Each polytope will be constructed from a polytope of dimension one lower. This is similar to the construction done by Bayer and Billera \cite{BB}, although, rather than bipyramids, we will build the dual, prisms.  For a $(d-1)$-polytope $Q$, we will denote a $d$-pyramid over it as $PQ$ and the d-prism over it as $B^{*}Q$, following Bayer and Billera's notation for pyramids and bipyramids but using $B^{*}Q$ to denote the dual of the bipyramid $BQ$. However, since we are interested in geometric aspects of the construction, we will fix the geometry of the polytopes and not just the combinatorics.  

The \emph{prism} $B^{*}Q$ is  $Q \times I$, where $I=[0,k]$ for some $k$. Then any $i$-face $F$ of $B^{*}Q$ is either an $i$-face of one of $Q\times \{0\}$ or $Q\times \{k\}$, or, for some $(i-1)$-face $G \subseteq Q$, $F=G\times I$, which is orthogonal to both $Q \times \{0\}$ and $Q\times \{k\}$. If $F$ is a face of this latter type, then $\alpha(F,B^{*}Q) = \alpha(G,Q)$.  No angles change as the distance between the two copies of $Q$ varies, so the angle sums do not depend on $k$. More specifically, we have the following equations on the $f$-vector and angle sums:
\begin{equation}\label{alpha and B*}
\begin{split}
f_{0}(B^{*}Q) & = 2f_{0}(Q), \\
f_{i}(B^{*}Q) & =  2f_{i}(Q) + f_{i-1}(Q) \text{  for } 1 \leq i \leq d,\\
\alpha_{i}(B^{*}Q)& =  \alpha_{i}(Q)+\alpha_{i-1}(Q) \text{  for } 0 \leq i \leq d.
\end{split}
\end{equation} 

\begin{figure}
\begin{center}
\includegraphics[width=5.0in]{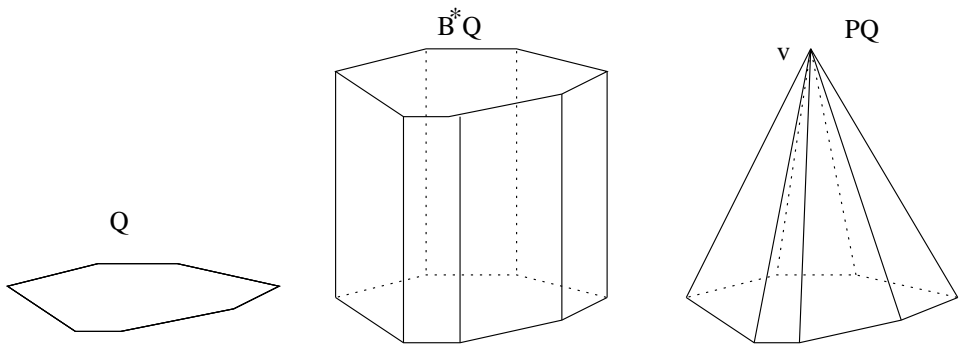}
\caption{The polytope $Q$; $B^{*}Q$, the prism over $Q$; and $PQ$, the pyramid over $Q$.}
\end{center}
\end{figure}
\vspace{.1in}

Now we define the \emph{pyramid} $PQ$.  We start by placing a $(d-1)$-dimensional polytope $Q$ in the hyperplane $x_{d}=0$ in $\R^{d}$. We then place a vertex $v$ along the line through the centroid of $Q$ and perpendicular to $Q$, so that it has $d$th coordinate $k>0$.  $PQ$ is then the convex hull of $v$ and $Q$.  An $i$-face of $PQ$ is either an $i$-face of $Q$ and therefore part of the base of the pyramid or the convex hull of $v$ and an $(i-1)$-face of $Q$.  We will refer to the latter as \emph{sides}. The angles formed between the sides and faces in the base increase as $k$ does. For this reason, we will denote the pyramid by $P_{k}Q$ to specify the height of $v$ and fix the geometry of the construction.  For any $k$, the pyramid operation has the following effect on the $f$-vector:
\begin{equation}\label{f and P}
f_{i}(PQ)  = f_{i}(Q) + f_{i-1}(Q), \text{  for } 0\leq i \leq d.
\end{equation}

\noindent We note two limiting cases of the pyramid operation: $P_{0}Q$, when $k$ tends to $0$ and $P_{\infty}Q$, when $k$ tends to infinity.  Although neither is actually a $d$-pyramid, one can easily find the limits of the angle sums as $k$ tends to $0$ or infinity, and we will define these values as the angle sums for $P_{0}Q$ and $P_{\infty}Q$.  Since the values of the angle sums vary continuously as $k$ does, we can find pyramids with angle sums that are arbitrarily close to those of $P_{0}Q$ and $P_{\infty}Q$.

For $P_{0}Q$, all angles made between the base and sides tend to $0$, so any interior angles at proper faces of
 the base are $0$. The interior angle at the base and at faces including $v$ are all $\frac{1}{2}$.
 Therefore, all the angles sums are dependent on the $f$-vector of the base $Q$. Then 
\begin{equation}\label{alpha P0}
\begin{split}
\alpha_i(P_{0}(Q)) & = \frac{1}{2} f_{i-1}(Q) \text{  for }0 \leq i \leq d-2,\\
\alpha_{d-1}(P_{0}(Q)) & = \frac{1}{2}f_{d-2}(Q) + \frac{1}{2}.
\end{split}
\end{equation}

For $P_{\infty}Q$, angles between the sides and base tend to right angles, so for any face $G \subseteq Q$, $\alpha(G,P_{\infty}Q)= \frac{1}{2} \alpha(G,Q)$.  For faces $F \subset P_{\infty}Q$ that are the convex hull of a face $G \subset Q$ and $v$, the interior angle at $F$ is the same as it was at $G$, that is, $\alpha(F,P_{\infty}Q) = \alpha(G, Q)$. Therefore
\begin{equation}\label{alpha and Pinfty}
\alpha_{i}(P_{\infty}Q)  =   \frac{1}{2} \alpha_{i}(Q)+\alpha_{i-1}(Q)
\text{  for }0 \leq i \leq d-1.
\end{equation}

We will sometimes want to iterate these constructions; we will write $C^{k}Q$ when we wish to apply a construction $C$ $k$ times in succession to $Q$. Taking a pyramid over a point $d$ times results in a $d$-simplex. Therefore, we will denote $d$-simplices as $P^d$, assuming a starting polytope of a point when one is not explicitly given.  

\section{The $\gamma$-vector}

In analogy to the $h$-vector, we define the \emph{$\gamma$-vector} as  \[\gamma(\P) = (\gamma_0(\P), \gamma_1(\P), \ldots, \gamma_d(\P)),\] where
\begin{equation}\label{gamma-vector}
\gamma_i(\P) = \sum_{j=0}^{i}(-1)^{i-j}{d-j \choose d-i}\alpha_{j-1}(\P).
\end{equation}
We note that $\gamma_0(\P) = 0$, $\gamma_1(\P) = \alpha_0(\P)$, and $\gamma_d(\P) = 1$ for all polytopes $\P$ and the transformation from the $\alpha$-vector to the $\gamma$-vector is invertible. 
We define

\hspace{.5in}$\gamma$-$f$-vector:
$(\gamma_0(P), \ldots, \gamma_d(P)| f_0(P) \ldots, f_d(P)),$ and

\hspace{.5in} $\gamma$-$h$-vector: $(\gamma_0(P), \ldots, \gamma_d(P)| h_0(P) \ldots, h_d(P)),$

\noindent and note that they are invertible linear transformations of the $\alpha$-$f$-vector. Therefore, a set of $\alpha$- or $\alpha$-$f$-vectors is affinely (linearly) independent if and only if the corresponding set of $\gamma$- or $\gamma$-$f$-vectors is.  Since the transformation between the $f$- and $h$-vector is also invertible and independent of the $\alpha$-vector, the affine (linear) independence of a set of $\alpha$-$f$-vectors is equivalent to the affine (linear) independence of the corresponding set of $\gamma$-$h$-vectors.

Kleinschmidt and Smilansky \cite{KS} defined a vector $\sigma_i(\D)$ that agrees with the $\gamma$-vector on spherical simplices.  The sphere was decomposed into regions by the great spheres that defined $\D$, and $\sigma_i(\D)$ measured the normalized area of all the regions that were reached from $\D$ by crossing $i$ great spheres.  We have chosen a different name for our vector to avoid confusion in the definition. 
 
As with the $h$-vector formulation of the Dehn-Sommerville relations, we can rewrite the Perles relations in terms of the $\gamma$-vector:
\begin{theorem}\label{h-Perles}
For a simplicial $d$-polytope $\P$, 
\[\gamma_i(\P) + \gamma_{d-i}(\P) = h_i(\P) \quad \text{for } 0 \leq i \leq d.\]
\end{theorem}
\begin{proof}
The proof follows the one given for Corollary 2.2 in \cite{BB}. 
We take the linear combination
\[\sum_{i=0}^{r}(-1)^i{d-i \choose d-r}S_{i-1}^d,\]
where $S_k^{d}$ is
\[\sum_{j=k}^{d-1} (-1)^{j} {j+1\choose k+1} \alpha_{j}(\P) =
(-1)^{d}(\alpha_{k}(\P) - f_{k}(\P)),\]
the $k$th Perles relation on simplicial $d$-polytopes.  On the right hand side, the sum becomes
\begin{equation*}
\begin{split}
\sum_{i=0}^r(-1)^i&{d-i \choose d-r}(-1)^{d}(\alpha_{i-1}(\P) - f_{i-1}(\P))\\
& = (-1)^{d-r}\left[\sum_{i=0}^r(-1)^{i-r}{d-i \choose d-r}\alpha_{i-1}(\P) - \sum_{i=0}^r(-1)^{i-r}{d-i \choose d-r}f_{i-1}(\P) \right]\\
& = (-1)^{d-r}\left(\gamma_r(\P) - h_r(\P)\right).
\end{split}
\end{equation*}
The left hand side is 
\begin{equation*}
\begin{split}
\sum_{i=0}^r(-1)^{i}{d-i \choose d-r}\sum_{j=i-1}^{d-1} (-1)^{j} &{j+1\choose i} \alpha_{j}(\P) \\
&= \sum_{i=0}^r(-1)^{i}{d-i \choose d-r}\sum_{j=i}^{d} (-1)^{j-1} {j\choose i} \alpha_{j-1}(\P)\\
&=\sum_{j=0}^{d} (-1)^{j-1}\alpha_{j-1}(\P)\sum_{i=0}^{j}(-1)^{i}{d-i \choose d-r} {j\choose i}.
\end{split}
\end{equation*}
Then we apply the identity 
\[\sum_{s=0}^n (-1)^{s}{s+m \choose t}{n \choose s} = (-1)^{n}{m \choose t-n}\] 
to simplify the interior sum
\begin{equation*}
\sum_{i=0}^{j}(-1)^{i}{d-i \choose d-r}{j\choose i}   = \sum_{s=0}^{j}(-1)^{j-s}{d-j+s \choose d-r} {j \choose s}= {d-j \choose r}.
\end{equation*}
Therefore the left hand side simplifies to 
\begin{equation*}
\begin{split}
\sum_{i=0}^r(-1)^{i}{d-i \choose d-r}\sum_{j=i-1}^{d-1} (-1)^{j} {j+1\choose i} \alpha_{j}(\P) & = \sum_{j=0}^{d} (-1)^{j-1}{d-j \choose r}\alpha_{j-1}(\P)\\
& = (-1)^{d-r+1}\gamma_{d-r}(\P).
\end{split}
\end{equation*}
Putting these results together we see that
\[\gamma_{d-r}(\P) =  h_r(\P)-\gamma_r(\P).\]
\end{proof}

In preparation for using the constructions to create affinely independent $\gamma$-vectors, we consider the effect of the pyramid and prism constructions on the $\gamma$-vector.  In the following proposition, we consider the $h$-vector entries strictly as a linear combination of the $f$-vector entries and do not assume that the polytope is simplicial.  
\begin{proposition}\label{gamma and P0}
If $Q$ is a $(d-1)$-polytope, 
\[h(PQ) =(h(Q), 1)\]
and
\[\gamma_i(P_{0}Q) = \left(0, \frac{1}{2}h_0,\frac{1}{2}h_1, \ldots, \frac{1}{2}h_{d-2},1\right).\]
\end{proposition}
\begin{proof}
The first relation is Proposition 3.1 of \cite{BB}.

For any polytope $Q$, $\gamma_0(Q) = 0$ and $\gamma_d(Q) = 1$.  From \eqref{alpha P0}, if Q is a $(d-1)$-polytope, $\alpha_j(P_{0}Q) = \frac{1}{2}f_{j-1}(Q)$ for $0 \leq j \leq d-2$.  Therefore, for $1 \leq i \leq d-1$,
\begin{equation*}
\begin{split}
\gamma_i(P_{0}Q) &= \sum_{j=0}^{i}(-1)^{i-j}{d-j \choose d-i}\alpha_{j-1}(P_{0}Q)\\
& = \frac{1}{2}\sum_{j=0}^{i-1}(-1)^{i-j-1}{d-j-1 \choose d-i}f_{j-1}(Q)\\
& = \frac{1}{2}h_{i-1}(Q).
\end{split}
\end{equation*}
\end{proof}
Induction using the proposition shows that that the $h$-vector of a $d$-simplex $\D$ is $(1,1,\ldots,1)$, so that $\gamma_i(\D) + \gamma_{d-i}(\D) = 1$ for $0 \leq i \leq d$ by Theorem \ref{h-Perles}.  In particular, the proposition shows that $\gamma(P_0^{d-1}P) = (0,\frac{1}{2}, \ldots \frac{1}{2}, 1)$.

\begin{proposition}\label{gamma and Pinfinity}
If $Q$ is a $(d-1)$-polytope, 
\[\gamma(P_{\infty}Q) = \frac{1}{2}\left[(0, \gamma(Q))+ (\gamma(Q),1)\right].\]
More generally,  
\[\gamma_i((P_{\infty})^{k}Q) = \frac{1}{2^k}\sum_{j=0}^{k}{k \choose j}\gamma_{i-j}(Q) \quad \text{for } 0 \leq i \leq d,\]
where
\[\gamma_k(Q) = \begin{cases}
1  & \text{ if }k \geq d\\
0 & \text{ if }k \leq 0.
\end{cases}\]
\end{proposition}
\begin{proof}
As for any $d$-polytope, $\gamma_0(P_{\infty}Q) = 0$ and $\gamma_d(P_{\infty}Q) = 1$. By \eqref{alpha and Pinfty}, $\alpha_{j}(P_{\infty}Q) =  \frac{1}{2} \alpha_{j}(Q)+\alpha_{j-1}(Q)$ for $0 \leq j \leq d-1$.
Then we calculate for $0 \leq i \leq d-1$:
\begin{equation*}
\begin{split}
\gamma_i(P_{\infty}Q) &= \sum_{j=0}^{i}(-1)^{i-j}{d-j \choose d-i}\alpha_{j-1}(P_{\infty}Q)\\ 
& = \frac{1}{2}\sum_{j=0}^{i}(-1)^{i-j}{d-j \choose d-i}\alpha_{j-1}(Q) \\
& \qquad \qquad \qquad + \sum_{j=0}^{i-1}(-1)^{i-j-1}{d-j-1 \choose d-i}\alpha_{j-1}(Q)\\
& = \frac{1}{2}\sum_{j=0}^{i-1}(-1)^{i-j-1}{d-j-1 \choose d-i}\alpha_{j-1}(Q) \\
& \qquad \qquad \qquad + \frac{1}{2}\sum_{j=0}^{i}(-1)^{i-j}{d-j-1 \choose d-i-1}\alpha_{j-1}(Q) \\
& = \frac{1}{2}\gamma_{i}(Q) + \frac{1}{2}\gamma_{i-1}(Q).
\end{split}
\end{equation*}

As we iterate the $P_{\infty}$ construction, $\gamma_i((P_{\infty})^kQ)$ is a linear combination of the set $\{\gamma_j(Q): i-k \leq j \leq i\}$ and the coefficient of $\gamma_j(Q)$ is half the sum of the coefficients of $\gamma_j(Q)$ and $\gamma_{j-1}(Q)$ in the linear combination for $\gamma_i((P_{\infty})^{k-1}Q)$.  Therefore, the coefficient of $\gamma_j(Q)$ in  $\gamma_i((P_{\infty})^kQ)$ is $\frac{1}{2^k}{k \choose j}$.
\end{proof}

\begin{proposition}\label{gamma and B*}
If $Q$ is a $(d-1)$-polytope, $\gamma(B^{*}Q) = (\gamma(Q), 1)$.
\end{proposition}
\begin{proof}
As for any $d$-polytope, $\gamma_d(B^{*}Q) = 1$.  By \eqref{alpha and B*}, $\alpha_{j}(B^{*}Q) =  \alpha_{j}(Q)+\alpha_{j-1}(Q)$ for $0 \leq j \leq d$.

Then we can calculate for $0 \leq i \leq d$:
\begin{equation*}
\begin{split}
\gamma_i(B^{*}Q) &= \sum_{j=0}^{i}(-1)^{i-j}{d-j \choose d-i}\alpha_{j-1}(B^{*}Q)\\ 
& = \sum_{j=0}^{i}(-1)^{i-j}{d-j \choose d-i}\alpha_{j-1}(Q) + \sum_{j=0}^{i-1}(-1)^{i-j-1}{d-j-1 \choose d-i}\alpha_{j-1}(Q)\\
& = \sum_{j=0}^{i}(-1)^{i-j}{d-j-1 \choose d-i-1}\alpha_{j-1}(Q)\\
& = \gamma_i(Q).
\end{split}
\end{equation*}
\end{proof}

\section{Spans of $\alpha$ and $\alpha$-$f$-vectors}
Using the $B^{*}$, $P_0$ and $P_{\infty}$ constructions, we will build families of polytopes with affinely independent
$\alpha$-vectors or $\alpha$-$f$-vectors.  We will use these families to span the spaces of $\alpha$-vectors and $\alpha$-$f$-vectors
defined by the Gram and Perles relations.  However, $P_0$ and $P_{\infty}$ are limiting cases of the pyramid construction and do not create $d$-polytopes.  The following lemma to guarantees a set of polytopes whose $\alpha$-vectors or $\alpha$-$f$-vectors have the same independence properties as those made by the constructions.
 
\begin{lemma}\label{backing off limiting}
Let $\e>0$  and $\{Q_i: i=0,\ldots, k\}$ be given,  where each $Q_i$ is a $d$-polytope or has form  $(P_0)^{k}Q$ or $(P_{\infty})^{k}(P_0)^{l}Q$ for some nonnegative integers $k$ and $l$ and a polytope $Q$ such that $k + l + \dim(Q) = d$.  

Then there is a set of $d$-polytopes $\{Q'_i: i=0,\ldots, k\}$, where $Q'_i$ and $Q_i$ have the same $f$-vector and $\left|\alpha_j(Q'_i)- \alpha_j(Q_i)\right|< \e$ for all $i$ and $j$.  Further, if the $\alpha$-vectors or $\alpha$-$f$-vectors of the $Q_i$ are affinely independent, $\e$ can be chosen so that the $\alpha$-vectors or $\alpha$-$f$-vectors of the $Q'_i$ are also affinely independent.
\end{lemma}
\begin{proof}
For each $i$ we will define constants $M_i$ and $\delta_i$.  Suppose $Q_i$ has form $P_{\infty}Q$.  Since the angle sums are continuous, we can choose $M_{i}^{j}$ so that $\alpha_j(P_{N}Q)$ is within $\e$ of $\alpha_j(P_{\infty}Q)$, for any $N \geq M_i^{j}$.  Let $M_{i}=\max_{j}M^{j}_{i}$.  Then for any $N \geq M_{i}$, $\left|\alpha_{j}(P_{N}Q) - \alpha_{j}(P_{\infty}Q)\right| < \e$ for all $j$. 

If $Q_i = (P_{\infty})^kQ$, we can iterate this process with difference $\e/k$.  Starting with $M_{i_0} = 1$, we iteratively choose $M_{i_{m}}\geq M_{i_{m-1}}$ by the same process as above so that  $\left|\alpha_{j}((P_{\infty})^{k-m}(P_{N})^m Q) - \alpha_{j}((P_{\infty})^{k-m+1}(P_{N})^{m-1} Q)\right| < \e/k$ for $N \geq M_{i_{m}}$ and all $j$.  Then we let $M_{i} = M_{i_k}$ so that if $N \geq M_i$,  $\left|\alpha_j((P_{N})^{k} Q) - \alpha_j((P_{\infty})^{k} Q)\right|<\e$.   

If $Q_i$ has form $P_{0}Q$, an analogous argument finds $\delta_i$ such that for all $\delta\leq \delta_i$ and $j$, $\left|\alpha_j(P_{\delta} Q) - \alpha_j(P_{0} Q)\right|<\e$.  Since the $\alpha$-vector of $P_{0}Q$ is entirely determined by the combinatorics of $Q$, this one step is also sufficient to choose $\delta_i$ for $Q_i = (P_{0})^{k}Q$.

Suppose $Q_i$ has form $(P_{\infty})^{k}(P_0)^{l}Q$ with $l \geq 1$.  First choose $\delta_i$ so that for all $\delta\leq \delta_i$, $\left|\alpha_j((P_{\delta})^{l} Q) - \alpha_j((P_{0})^{l} Q)\right|<\frac{\e}{2^{k+1}}$ for all $j$.  Then by \eqref{alpha and Pinfty} the $P_{\infty}$ construction will less than double any differences in angle sum values, so \[\left|\alpha_j((P_{\infty})^k(P_{\delta})^{l} Q) - \alpha_j((P_{\infty})^k(P_{0})^{l} Q)\right|<\e/2\]
for $\delta\leq \delta_i$ and all $i$ and $j$. Next choose $M_i$ so that  
\[\left|\alpha_{j}((P_{N})^k(P_{\delta})^{l}Q) - \alpha_{j}((P_{\infty})^k(P_{\delta})^{l}Q)\right| < \e/2\] for $N \geq M_{i}$ and all $i$ and $j$.  Then \[\left|\alpha_{j}((P_{N})^k(P_{\delta})^{l}Q) - \alpha_{j}((P_{\infty})^k(P_{0})^{l}Q)\right| < \e\]
 for $N \geq M_{i}$, $\delta\leq \delta_i$,  and all $j$.  

Now we choose
\[Q'_i = \begin{cases}
(P_{\delta_i})^{k}Q & \text{ if } Q_i = (P_{0})^{k}Q\\
(P_{M_i})^k(P_{\delta_i})^{l}Q & \text{ if } Q_i = (P_{\infty})^{k}(P_{0})^{l}Q\\
Q_i & \text{ if } Q \text{ is a }d\text{-polytope}.
\end{cases}\]
The $f$-vectors of $Q_i$ and $Q'_i$ are the same since they are pyramids of the same degree over the same polytope and $\left|\alpha_j(Q'_i)- \alpha_k(Q_i)\right|< \e$ for all $i$ and $j$. 

Since affine independence is an open condition, if the $Q_i$ have affinely independent $\alpha$-vectors or $\alpha$-$f$-vectors, we can choose $\e$ small enough that the $Q'_i$ given above have affinely independent $\alpha$-vectors or $\alpha$-$f$-vectors, respectively.  
\end{proof}

The $P_{\infty}$ construction will be useful for increasing the dimension of a set of polytopes and maintaining the affine independence of their $\alpha$-vectors.
\begin{lemma}\label{P_{infty} independence}
If a set of $(d-1)$-polytopes $\{Q_i: i = 0,\ldots, k\}$ has affinely independent $\alpha$-vectors, then the $\alpha$-vectors of the set $\{P_{\infty}Q_i: i = 0,\ldots, k\}$ are also affinely independent.
\end{lemma}
\begin{proof}
We will work with the $\gamma$-vectors for ease of computation. Since the last entry of $\gamma(Q)$ is 1 for every polytope, the affine independence of a set of $\gamma$-vectors is equivalent to their linear independence.   

Based on Proposition \ref{gamma and Pinfinity}, we can write
\begin{equation*}
\gamma(P_{\infty}Q) =\textbf{A}\begin{bmatrix} \gamma(Q)\\1\end{bmatrix}
\end{equation*}
where
\begin{equation}\label{Pinfty trans on gamma}
\textbf{A} = \frac{1}{2}\left[\begin{array}{ccccc}
 1& & & &  \\
 1&1& & &  \\
  &1&1& &   \\
 &&\ddots&\ddots& \\ 
  & & &1&1\\
\end{array}\right],
\end{equation}
a $(d+1)\times(d+1)$ matrix where all other entries are 0.

Clearly $\textbf{A}$ is invertible.  Its inverse is
\begin{equation}\label{Ainverse}
\textbf{A}^{-1} = \left[\begin{array}{rrrrr}
 2 & & & &\\ 
 -2 & 2 &  & &   \\
 2 & -2 & 2 &&\\
\vdots& & &\ddots\\ 
&\ldots&2&-2&2\\
\end{array}\right],
\end{equation}
where all entries on and below the diagonal alternate between $2$ and $-2$ and entries above the diagonal are 0. 
Since the matrix is invertible, the $P_{\infty}$ transformation preserves the linear independence of the $\gamma$-vectors of $Q_{i}$ for $i=1,\ldots,k$.
\end{proof}

\begin{theorem}\label{Angle sums of simplices} The affine span of the $\alpha$-vectors of $d$-simplices is
$\left\lfloor\frac{d-1}{2}\right\rfloor$-dimensional.
\end{theorem}
\begin{proof}   As in Lemma \ref{P_{infty} independence}, we will work with the $\gamma$-vector.  Let $A$ be the affine space spanned by the $\gamma$-vectors of $d$-simplices. 

If $\D$ is a $d$-simplex, $\gamma_k(\D) + \gamma_{d-k}(\D) = 1.$ These relations are clearly independent for $i=0,\ldots,\left\lfloor\frac{d}{2}\right\rfloor$.  Since the $\gamma$-vector is $(d+1)$-dimensional and all $\alpha$-vectors lie in the plane $\gamma_{d} = 1$, 
\begin{equation*}
\dim(A) \leq d+1-\left(\left\lfloor\frac{d}{2}\right\rfloor+1\right)- 1 =\left\lfloor\frac{d-1}{2}\right\rfloor.
\end{equation*}

We will prove that $\dim(A) \geq \left\lfloor\frac{d-1}{2}\right\rfloor$ by constructing a set of $\left\lfloor\frac{d-1}{2}\right\rfloor + 1$ simplices whose $\gamma$-vectors are affinely independent. 
  The proof will proceed by induction on $d$.  For $d=1$ and $d=2$, a line segment and a triangle (denoted $P$ and $P^{2}$, respectively) provide the one element needed for the basis. 

Let $d \geq 3$. Suppose we have a set of $\left\lfloor \frac{d-3}{2} \right\rfloor +1 = \left\lfloor \frac{d-1}{2} \right\rfloor$ $(d-2)$-simplices $\{Q_i: i = 1,\ldots,\left\lfloor \frac{d-1}{2} \right\rfloor\}$ with linearly independent $\gamma$-vectors. We claim that the vectors $\gamma\Big(P_0^{d-1}P\Big), \gamma\Big(P_{\infty}^2Q_1\Big), \ldots, \gamma\Big(P_{\infty}^2Q_{\left\lfloor \frac{d-1}{2} \right\rfloor}\Big)$ are linearly independent.  

We know the vectors $\gamma\Big(P_{\infty}^2Q_1\Big), \ldots, \gamma\Big(P_{\infty}^2Q_{\left\lfloor \frac{d-1}{2} \right\rfloor}\Big)$ are linearly independent by applying Lemma \ref{P_{infty} independence} twice.  We will show that adding the vector $\gamma(P_0^{d+1}P)$ increases the linear span by showing that the inverse image of $\gamma(P_0^{d+1}P)$ under the $P_{\infty}^2$ transformation is not in the linear span of the $\gamma(Q_i)$. 

Now by Proposition \ref{gamma and P0}, $\gamma(P_0^{d-1}P) = \left(0, \frac{1}{2},\frac{1}{2}, \ldots, \frac{1}{2}, 1\right)$.  Then using $\textbf{A}^{-1}$ from \eqref{Ainverse}, 
\[\left(\textbf{A}^{-1}\right)^2\left(\gamma(P_0^{d-1}P)\right) = \begin{bmatrix} 0\\2\\ \vdots\end{bmatrix},\]
where the later entries alternate in sign. Therefore, the last three entries do not have the same value. But each of the vectors \[\left(\textbf{A}^{-1}\right)^2\left(\gamma\left(P_{\infty}^2Q_i\right)\right) = \begin{bmatrix} \gamma(Q_i)\\1\\1\end{bmatrix}\]
has last three entries 1, so any vector in the span of the $\gamma(Q_i)$ has the same value on the last three entries.  Therefore $\left(\textbf{A}^{-1}\right)^2\left(\gamma(P_0^{d-1}P)\right)$ is outside the linear span of the $\gamma(Q_i)$ for $i=1, \ldots, \left\lfloor \frac{d-1}{2} \right\rfloor$.

This gives us a set of $\left\lfloor \frac{d-1}{2} \right\rfloor + 1$ $d$-simplices constructed by $P_{0}$ and $P_{\infty}$ operations whose $\alpha$-vectors are affinely independent. By Lemma \ref{backing off limiting}, there is a set of $\left\lfloor \frac{d-1}{2} \right\rfloor + 1$ $d$-simplices with affinely independent $\alpha$-vectors. 
\end{proof}
 
Using what is known about the affine span of the $f$-vectors of simplicial polytopes together with the results of the preceding theorem, we can determine the affine span of the $\alpha$-$f$-vectors of simplicial polytopes. It is appropriate to consider this vector rather than the $\alpha$-vector in describing the angles of simplicial polytopes since the Perles relations refer to face numbers as well as angle sums.
\begin{theorem}\label{Angle sums of simplicial polytopes} The affine span of the $\alpha$-$f$-vectors of simplicial $d$-polytopes has dimension $d-1$.  The space is spanned by $\left\lfloor\frac{d+1}{2}\right\rfloor$ simplices, as in Theorem \ref{Angle sums of simplices}, and $\left\lfloor\frac{d}{2}\right\rfloor$ non-simplices which are combinatorially independent simplicial polytopes. 
\end{theorem}
\begin{proof} We will work with the $\gamma$-$h$-vectors.  Let $A_S$ be the affine space spanned by the $\gamma$-$h$-vectors of simplicial polytopes.  

By the argument in Theorem \ref{Angle sums of simplices}, there are $\lfloor \frac{d}{2}\rfloor+1$ Perles relations that are independent with regard to angle sums, and in this case each includes a different element of the $h$-vector.  The other relation on the angle sums is that $\gamma_d(P)=1$ for all polytopes.  Similarly, the Dehn-Sommerville relations on the $h$-vector show that there are $\lfloor \frac{d+1}{2}\rfloor$ independent Dehn-Sommerville relations.  Since these relations include no angle sums, they are independent of the Perles relations.  Therefore, 
\begin{equation*}
\text{dim}(A_S) \leq 2d+2-\left(\left\lfloor\frac{d}{2}\right\rfloor+1\right)- 1 - 
\left(\left\lfloor \frac{d+1}{2}\right\rfloor\right) - 1=d-1.
\end{equation*}

The affine span of the $h$-vectors of simplicial $d$-polytopes has dimension $\left\lfloor\frac{d}{2}\right\rfloor$.  In Bayer and Billera \cite{BB}, a set of $\left\lfloor\frac{d}{2}\right\rfloor + 1$ simplicial polytopes with affinely independent $h$-vectors is given, spanning the space defined by the Dehn-Sommerville equations.  This basis includes one simplex.  We can combine the $\left\lfloor\frac{d}{2}\right\rfloor$ non-simplices of this basis with the $\left\lfloor\frac{d+1}{2}\right\rfloor$ simplices given in Theorem \ref{Angle sums of simplices}.  

If the $\gamma$-$h$-vectors of this set are affinely dependent, the dependency must occur within the set of simplices since otherwise this would give a dependency on the $h$-vectors of the set of polytopes in \cite{BB}.  However, this would be a contradiction to the previous theorem.  Therefore the $\gamma$-$h$-vectors of the set of polytopes are affinely independent.
\end{proof} 

For the $\alpha$-$f$-vectors of simplicial polytopes, this shows that the dimensions beyond those determined combinatorially are found in variation of the angle sums of simplices.  This means that once we have considered the degrees of freedom in the geometry of the simplex, all other degrees of freedom in the geometry of simplicial polytopes can be described by variation in combinatorial dimensions.

We can similarly build a set of polytopes whose $\alpha$-$f$-vectors span the space defined by the Gram and Euler relations. We will use a method similar to the proof of Theorem \ref{Angle sums of simplices}, but first we prove the following lemma.
\begin{lemma}\label{inverse image of P on B*Q}
Let $Q$ be a $d$-polytope with $f$-vector $f = (f_0, f_1, \ldots, f_d)$ and let $\f = (1, \f_0, \f_1, \ldots, \f_{d-1})$ be the inverse image of $f$ under the pyramid transformation.  Also, let the $(d+1)$-polytope $B^{*}Q$ have $f$-vector $f^{*} = (f^{*}_0, f^{*}_1, \ldots, f^{*}_{d+1})$ and inverse image $\f^{*} = (1, \f^{*}_0, \f^{*}_1, \ldots, \f^{*}_{d})$ under the pyramid transformation.  Then \[\sum_{i=0}^{d} (-1)^i\f^{*}_i = \sum_{i=0}^{d-1} (-1)^i\f_i +1.\]
\end{lemma}
\begin{proof}
By \eqref{f and P}, if we extend the $f$-vector of a polytope $Q$ to $(1, f(Q))$,
\begin{equation}\label{Matrix B}
f(PQ) = \textbf{B}\begin{bmatrix} 1\\f(Q)\end{bmatrix},
\end{equation}
where
\begin{equation*}
\textbf{B} = \left[\begin{array}{ccccc}
 1&1&&&\\
 &1&1&&\\
 &&\ddots&\ddots&\\
 &&&1&1\\
 &&&&1
\end{array}\right],
\end{equation*}
a $(d+1)\times(d+1)$ matrix where all other entries are 0.  Therefore the matrix for the inverse transformation is
\begin{equation}\label{Binverse}
\textbf{B}^{-1} = \left[\begin{array}{rrrrr}
 1&-1&1&-1&\ldots\\
 &1 & -1 & 1 & \ldots\\
 &&\ddots&\ddots&\\
 &&&1&-1\\
 &&&&1
\end{array}\right],
\end{equation}
where all the entries below the main diagonal are 0.  Multiplication by $\textbf{B}^{-1}$ shows 
\[\f_i = \sum_{j=i+1}^d (-1)^{i-j+1} f_j \quad \text{for } 0 \leq i \leq d-1\]
 and
 \begin{equation}\label{barf* to f*}
 \f^{*}_i = \sum_{j=i+1}^{d+1} (-1)^{i-j+1} f^{*}_j \quad \text{for } 0 \leq i \leq d.
\end{equation}
Since $f$ and $f^{*}$ are $f$-vectors of polytopes, $f_d = f^{*}_{d+1} = 1$.  We also know from \eqref{alpha and B*} that $f^{*}_i = 2f_i + f_{i-1}$ for $1 \leq i \leq d$.  Therefore, we can rewrite \eqref{barf* to f*} for $0 \leq i \leq d$ as:
\begin{equation*}
\begin{split}
\f^{*}_i &= 2\sum_{j=i+1}^{d} (-1)^{i-j+1}f_j + \sum_{j=i}^{d-1} (-1)^{i-j}f_{j}+(-1)^{i-d}f^{*}_{d+1}\\
& =f_i + \sum_{j=i+1}^{d} (-1)^{i-j+1}f_j \\
& = f_i + \f_i.
\end{split}
\end{equation*}
Now, taking the alternating sum we get
\begin{equation*}
\begin{split}
\sum_{i=0}^{d} (-1)^i\f^{*}_i & =\sum_{i=0}^{d} (-1)^i f_i + \sum_{i=0}^{d} (-1)^i \f_i\\
& = 1 + \sum_{i=0}^{d-1} (-1)^i \f_i,
\end{split}
\end{equation*}
where the last equality follows by the Euler relation on $Q$.
\end{proof}

\begin{theorem}\label{Angle sums of general polytopes} The affine span of the $\alpha$-$f$-vectors of general $d$-polytopes has dimension $2d-3$ for $d \geq 2$.
\end{theorem}
\begin{proof}
We will work with the set of $\gamma$-$f$-vectors. The Euler and Gram equations provide two independent equations on the $\gamma$-$f$-vectors.  We also know that $\gamma_{0}(P)=0$, $\gamma_{d}(P) = 1$,  and $f_{d}(P)=1$ for all polytopes $P$.  As long as $d>1$, these equations are independent.  Therefore, the span of the $\gamma$-$f$-vectors is at most $2d+2-5=2d-3$ if $d\geq 2$. To show this whole space is spanned, we will again proceed inductively on $d$.

The statement is true in two dimensions, since the $\gamma$-$f$-vectors of the triangle and the square (denoted $P^{2}$ and $B^{*}P$, respectively) are $\left(0,\frac{1}{2}, 1,3,3,1\right)$ and $(0,1,1,4,4,1)$. 

Suppose the statement is true for dimension $d-1$.  That is, there are $2(d-1)-2=2d-4$ affinely independent $\gamma$-$f$-vectors of $(d-1)$-polytopes: $Q_{1}, Q_{2}, \ldots, Q_{2d-4}$. Then we claim that the $\gamma$-$f$-vectors of the following polytopes are affinely independent:
\begin{equation}\label{list of general polytopes}
P_{\infty}Q_{1}, P_{\infty}Q_{2}, \ldots, P_{\infty}Q_{2d-4},(B^{*})^{d-2}P^{2} \text{ and }(B^{*})^{d-1}P
\end{equation} 

The affine independence of a set of $\gamma$-$f$-vectors is equivalent to their linear independence since $\gamma_d = f_d = 1$ for all $d$-polytopes.  Therefore, we will show that the $\gamma$-$f$-vectors of the polytopes in \eqref{list of general polytopes} are linearly independent. 

We will consider the effect of the $P_{\infty}$ construction on the $\gamma$-$f$-vector of a $(d-1)$-polytope $Q$.  We extend the $\gamma$-$f$-vector to the $(2d+2)$-vector $(\gamma(Q), 1, 1, f(Q))$, thinking of the additional entries as $\gamma_{d}(Q)$ and $f_{-1}(Q)$, respectively. Then we can write
\begin{equation}\label{Pinfty trans on gamma-f}
\gamma\text{-}f(P_{\infty}Q) = \textbf{C}\begin{bmatrix} \gamma(Q)\\1\\1\\f(Q)\end{bmatrix}
\end{equation}
where
\begin{equation*}
\textbf{C}=  \left[\begin{array}{c|c} \textbf{A}& 0 \\\hline 0 & \textbf{B}\end{array}\right],
\end{equation*}
a $(2d+2)\times(2d+2)$ matrix, with blocks $\textbf{A}$ \eqref{Pinfty trans on gamma} and $\textbf{B}$ \eqref{Matrix B}. Since this is an invertible matrix, the $\gamma$-$f$-vectors of $P_{\infty}Q_{i}$, $i=1,\ldots,2d-4$, are all linearly independent.  Therefore, we consider the $\gamma$-$f$-vectors of $(B^{*})^{d-2}P^2$ and $(B^{*})^{d-1}P$.

As in the proof of Theorem \ref{Angle sums of simplices}, we will consider the inverse images of the $\gamma$-$h$-vectors of these two polytopes in the $P_{\infty}$ transformation. Let 
\[v_1:=\textbf{C}^{-1}\left(\gamma\text{-}f((B^{*})^{d-2}P^2)\right)^{T} \qquad \text{and} \qquad
v_2:=\textbf{C}^{-1}\left(\gamma\text{-}f((B^{*})^{d-1}P)\right)^{T} \]
and consider these vectors in relation to the span of the $(\gamma(Q_i), 1, 1, f(Q_i))$.

The $f_{-1}$ entry of $\textbf{C}^{-1}\left(f(Q)\right)$ for any $d$-polytope $Q$ is 1 by the Euler relation, so the $f_{-1}$ entries of $v_1$ and $v_2$ are 1.  Similarly, for each $Q_i$, the alternating sum of the values of $f(Q_i)$ is 1 by the Euler relation and the $f_{-1}$ entry is 1 as well. Therefore, if $v_1$ and $v_2$ are in the linear span of the $\gamma$-$f(Q_i)$, the entries $\alpha_d$, $f_{-1}$, and the corresponding alternating sum must be 1.  

However, by Lemma \ref{inverse image of P on B*Q}, the alternating sum of the entries $f_0, f_1, \ldots, f_{d-1}$ of $\textbf{C}^{-1}\left(f\left(B^{*}Q\right)\right)$ for a $(d-1)$-polytope $Q$ is one greater than the alternating sum of the entries $f_0, f_1, \ldots, f_{d-2}$ of $\textbf{C}^{-1}\left(f\left(Q\right)\right)$. This alternating sum is 1 for $Q=P$ or $Q=P^2$, so the alternating sum of the entries $f_0, f_1, \ldots, f_{d-1}$ of $v_1$ and $v_2$ must be at least 2.  Therefore neither can be a linear combination of the extended $\gamma$-$f$-vectors of $Q_{i}$ for $i=1,\ldots,2d-4$. 

Now we consider the $\gamma$-vector entries.
By Proposition \ref{gamma and B*},  $\gamma((B^{*})^{d-2}P^2) = \left(0, \frac{1}{2}, 1, \ldots,1\right)$ and
$\gamma((B^{*})^{d-1}P) = \left(0,1,\ldots, 1\right).$
Therefore the $\gamma$-vector portions of $v_1$ and $v_2$ are
$(0,1, \ldots ,1)$ and $(0,2,0 ,2,\ldots)$, respectively.  Then $\gamma_d \neq f_{-1}$ in $v_2$, even though $\gamma_d =f_{-1}$  for $v_1$ and each of the vectors $(\gamma(Q_i), 1, 1, f(Q_i))$.  Therefore $v_2 \notin \spanr\{(\gamma(Q_i), 1, 1, f(Q_i)), v_1\}$
and the polytopes in \eqref{list of general polytopes} have affinely independent $\alpha$-$f$-vectors.

Then by Lemma \ref{backing off limiting}, we know that we have a set of $d$-polytopes of size $2d-2$ with affinely independent $\alpha$-vectors. 
\end{proof}

We note that the set of polytopes which span the space of $\alpha$-$f$-vectors has significant duplication in the $\alpha$-vectors.  For instance, the polytopes $P_{\infty}(B^{*})^{k}P$ and $(B^{*})^{k}P^{2}$ have the same angle sums for all $k \geq 1$.

These results strengthen the correspondence between the geometric structure and the combinatorial structure of polytopes.  The Gram and Perles relations are close analogs of the Euler and Dehn-Sommerville relations.  In this paper, we have shown that the affine dimensions closely correspond.  The affine span of the $\alpha$-vectors of $d$-simplices has the same dimension as the span of the $f$-vectors of simplicial $(d-1)$-polytopes.  Also, the affine span of the $\alpha$-$f$-vectors of simplicial $d$-polytopes has the same dimension as the span of the $f$-vectors of $d$-polytopes.  It would be interesting to speculate whether there is a deeper significance to this relationship. 

The use of the $\gamma$-vector also raises questions about the nature of this measure on angle sums.  For the $h$-vectors of simplicial polytopes, there are many results bounding the values.  The Upper Bound Theorem \cite{McM3} bounds the $h$-vector entries above by those of the cyclic polytope of the same dimension and the same number of vertices.  The Generalized Lower Bound Theorem \cite{McMWa} shows that the first $\left\lfloor \frac{d}{2} \right \rfloor$ entries of the $h$-vector are unimodal and the $g$-Theorem gives bounds on the differences between adjacent entries of the vector \cite{BL, Stanley}.  In the case of the $\gamma$-vector the bounds on its entries are unexplored.  Initial examples show that the $\gamma$-vector may be more tractable on non-simplicial polytopes than the $h$-vector; for example, the basis polytopes for the theorems in this chapter, many of which are not simplicial, all have monotonic $\gamma$-vectors.  This is not the case for all polytopes (for example, the bipyramid made by gluing two regular tetrahedra along a face), but unimodality may be true in general and monotonicity in specific cases such as for simplices.

\section*{Acknowledgements}
This work was done with the support and guidance of my advisor, Louis J. Billera.  I am also very thankful for the careful reading and advice of David W. Henderson.


\begin{thebibliography}{99}

\bibitem{BB} Bayer, M.M., Billera, L.J.; Counting faces and chains in polytopes and posets, in \emph{Combinatorics and Algebra}, Greene, C., (ed.), Contemporary Mathematics \textbf{34}. Providence: Amer. Math. Soc., 1984, 207-250.

\bibitem{BB2} Bayer, M.M., Billera, L.J.; Generalized
Dehn-Sommerville relations for polytopes, spheres and Eulerian partially ordered sets. \emph{Invent. Math.}, \textbf{79}(1), 1985, 143-157.

\bibitem{BL} Billera, L.J., Lee, C.W.; A proof of the sufficiency of McMullen's conditions for $f$-vectors of simplicial convex polytopes. \emph{J. Comb. Theory A}, \textbf{31}(1981), 237-255.

\bibitem{G}  Gram, J.P.; Om rumvinklerne i et polyeder. \emph{Tidsskrift for Math. (Copenhagen)}, \textbf{4}(3), 1874, 161-163.

\bibitem{Gr} Gr\"unbaum, B.; \emph{Convex Polytopes}, Graduate Texts in Mathematics \textbf{221}. Springer-Verlag: New York, 2003.

\bibitem{Hohn} H\"ohn, W.; \emph{Winkel und Winkelsumme im n-dimensionalen Euklidischen Simplex}.  Ph.D. Thesis. E.T.H. Zurich, 1953.

\bibitem{KS} Kleinschmidt, Peter, Smilansky, Zeev; New results for simplicial spherical polytopes, in \emph{Discrete and Computational Geometry (New Brunswick, NJ, 1989/1990)}, DIMACS Ser. Discrete Math. Theoret. Comput. Sci. \textbf{6}, Amer. Math. Soc., Providence, RI, 1991, 187-197.

\bibitem{McM3} McMullen, P.; The maximum numbers of faces of a convex polytope. \emph{Mathematika}, \textbf{17}(1970), 179-184.

\bibitem{McMWa} McMullen, P., Walkup, D.W.; A generalized lower bond conjecture for simplicial polytopes.  \emph{Mathematika}, \textbf{18}(1971), 264-273.

\bibitem{PS} Perles, M.A., Shephard, G.C.; Angle Sums of Convex Polytopes. \emph{Math. Scand.}, \textbf{21}(1967), 199-218.

\bibitem{S}  Shephard, G.C.; An Elementary Proof of Gram's Theorem for Convex Polytopes. \emph{Can. J. Math.}, \textbf{19}(1967), 1214-1217.

\bibitem{So}  Sommerville, D.M.Y.; The relations connecting the angle-sums and the volume of a polytope in space of n dimensions. \emph{Proc. Roy. Soc. London, Ser. A}, \textbf{115}(1927), 103-119.
\bibitem{Stanley} Stanley, R.; The number of faces of a simplicial convex polytope. \emph{Advances in Math.}, \textbf{35}(1980), 236-238.

\bibitem{W}  Welzl, E.; Gram's equation-a probabilistic proof, in \emph{Results and Trends in Theoretical Computer Science (Graz, 1994)}, Lecture Notes in Comput. Sci., \textbf{812}, . Berlin: Springer, 1994, 422-424.

\bibitem{Ziegler} Ziegler, G.M.; \emph{Lectures in Polytopes}, Graduate Texts in Mathematics \textbf{152}. Springer-Verlag: New York, 1995.
\end{thebibliography}
\end{document}